\documentclass{amsart}
\usepackage{amsmath,amssymb,amsfonts,latexsym,a4wide,amsthm}
\usepackage{amscd}  
\usepackage[all]{xy}

\title{A categorification of integral Specht modules}
\author{Mikhail Khovanov}
\address{M. K.: Department of Mathematics, Columbia University (USA)}
\email{khovanov\symbol{64}math.columbia.edu}
\thanks{M.K. was partially supported by the NSF grant DMS-0407784.} 
\author{Volodymyr Mazorchuk}
\address{V. M.: Department of Mathematics, Uppsala University 
(Sweden).}
\email{mazor\symbol{64}math.uu.se}
\thanks{V. M. was supported by STINT, the Royal Swedish Academy of 
Sciences, 
the
Swedish Research Council and the MPI in Bonn.}
\author{Catharina Stroppel}
\address{C. S.: Department of Mathematics, University of Glasgow 
(United 
Kingdom).}
\email{c.stroppel\symbol{64}maths.gla.ac.uk}
\thanks{C. S. was supported by EPSRC grant 32199}

\newtheorem{prop}{Proposition}

\newtheorem{theorem}[prop]{Theorem} 
\newtheorem{definition}{Definition} 

\newtheorem{Example}{Example}
\newtheorem{Remark}[prop]{Remark}
\newcommand{\oplusop}[1]{{\mathop{\oplus}\limits_{#1}}}

\begin{document}

\baselineskip 14pt

\def\la{\lambda}
\def\op{\operatorname}
\def\C{\mathbb C}
\def\R{\mathbb R}
\def\N{\mathbb N}
\def\Z{\mathbb Z}
\def\Q{\mathbb Q}
\def\g{\mathfrak g}
\def\p{\mathfrak p}
\def\h{\mathfrak h}
\def\n{\mathfrak n}
\newcommand{\mC}{\mathbb{C}}
\newcommand{\Ext}{\operatorname{Ext}}
\newcommand{\End}{\operatorname{End}}
\newcommand{\add}{\operatorname{add}}

\def\F{\mathbb F}
\def\S{\mathbb S}
\def\l{\lbrace}
\def\r{\rbrace}
\def\o{\otimes}
\def\lra{\longrightarrow}
\newcommand{\ba}{\mathbf{a}} 
\newcommand{\cA}{\mathcal{A}}
\newcommand{\cB}{\mathcal{B}}
\newcommand{\dmod}{\mathrm{-mod}}
\newcommand{\mc}{\mathcal}
\def\cF{\mathcal{F}}
\def\Hom{\textrm{Hom}}
\def\drawing#1{\begin{center} \epsfig{file=#1} \end{center}}
\def\mc{\mathcal}
\def\mf{\mathfrak}
\def\mb{\mathbb}
 
\def\yesnocases#1#2#3#4{\left\{ \begin{array}{ll} #1 & #2 \\ #3 & #4
\end{array} \right. }
 
\newcommand{\define}{\stackrel{\mbox{\scriptsize{def}}}{=}}
\def\hsm{\hspace{0.05in}}
 
\def\cO{\mc{O}}   
\def\cC{\mc{C}}
\def\sln{\mathfrak{sl}(n)}

\begin{abstract}
We suggest a simple definition for categorification of modules over rings 
and illustrate it by categorifying integral Specht modules over
the symmetric group and its Hecke algebra via the action of translation 
functors on some subcategories of category $\mathcal{O}$
for the Lie algebra $\mathfrak{sl}_n(\mathbb{C})$.  
\end{abstract}

\maketitle


\section{Introduction}
In this paper we describe categorifications of irreducible modules for the 
symmetric group $S_n$ and the corresponding (generic) 
Hecke algebra $\mathcal{H}_n$. The category of complex $S_n$-modules is 
semisimple and the 
irreducible modules are the so-called Specht modules ${\mathcal{S}}(\lambda)$,
 indexed by partitions $\la$ of $n$. Analogous statements hold for the Hecke
 algebra $\mathcal{H}_n$, where we denote the Specht module by
 $\tilde{\mathcal{S}}(\lambda)$.  The symmetric group $S_n$ is isomorphic to 
the Weyl group of the semi-simple complex Lie algebra $\mathfrak{sl}_n$. 
Associated to any composition $\mu$ of $n$, there is the corresponding Young 
subgroup $S_\mu$ of $S_n$ and the corresponding parabolic subalgebra $\p_\mu$ 
of $\mathfrak{sl}_n$. There is also the partition $\la$ associated with $\mu$
 (i.e. $\la$ coincides with $\mu$ up to reordering of the parts), and the 
dual partition $\la'$ of $\la$ with the corresponding Specht module
 ${\mathcal{S}}(\lambda')$.

On the other hand, there is the BGG category $\cO$ for $\mathfrak{sl}_n$  and its direct 
summand $\cO_0$ determined by the trivial central character. 
Let $\cO_0^\mu$ be the full subcategory of  $\cO_0$ consisting of 
 all $\p_\mu$-locally finite modules. 

We establish  a canonical bijection between the isomorphism 
classes of indecomposable pro\-jec\-tive-\-in\-jec\-tive 
modules in $\cO_0^\mu$
and the Kazhdan-Lusztig basis of ${\mathcal{S}}(\lambda')$,
 defined in \cite{KL}. 
The category  $\cO_0^\mu$ is equivalent to the category of finitely 
generated $A$-modules, where $A$ is a certain finite-dimensional 
positively graded Koszul algebra (\cite{BGS}). By working with  
the category $\cA$ of {\it graded} finitely-generated $A$-modules 
we obtain a canonical bijection 
between the isomorphism classes of indecomposable projective-injective 
modules in $\cA,$ up to grading shifts, and the Kazhdan-Lusztig basis of 
$\tilde{\mathcal{S}}(\lambda').$  

Every Specht module $\tilde{\mathcal{S}}(\lambda)$ has a symmetric, 
non-degenerate, $\mathcal{H}_n$-invariant bilinear form $<,>$ with values 
in $\Z[v,v^{-1}]$, which is unique up to a scalar (\cite{Murphy}). Under 
the bijection above, this form can be {\it categorified} and becomes 
the bifunctor 
$\operatorname{Hom}_A(_-, \operatorname{d}(_-))$ 
which takes two graded $A$-modules $M$, $N,$ with $M$ projective, 
to the 
graded vector space $\operatorname{Hom}_A(M,\operatorname{d}(N)),$ 
where $d$ is the standard duality functor on $\cA.$ The bilinear 
form, evaluated at images of $M$ and $N$ in the Grothendieck group, 
 equals the Hilbert polynomial of this graded vector space.  

To make sense of the  $\mathcal{H}_n$-invariance, we consider an action 
of the Hecke algebra on the category $\cA$ via the so-called translation 
functors in $\cO$, introduced by Jantzen (\cite{Ja1},\cite{BG}). 

As a result we get a {\it categorification} of the Specht modules for both
 the symmetric group (via the category $A$-mod) and the Hecke algebra (via 
the category $\cA$), where the action of
 the group ring and the Hecke algebra respectively is given by translation 
functors and the bilinear form is obtained from the bifunctor $\op{Hom}$ as 
briefly explained above. In particular, dimensions of endomorphism rings of 
projective-injective modules in $\cO_0^\mu$ can be obtained purely 
combinatorially using the non-degenerate form. 

We start by giving a definition of a categorification and then prove 
the two categorification results with the bilinear form obtained at the end. 
By a "module" we always mean a right module.

\section{Categorification of modules with integral structure constants} 

The Grothendieck group $K_0(\mc{B})$ of an abelian category $\mc{B}$ 
is generated by symbols $[M]$, where $M\in Ob(\mc{B})$, 
subject to the defining relations  $[M_2]=[M_1]+ [M_3]$ whenever there is a short exact sequence $0 \lra M_ 1\lra M_2 \lra M_3 \lra 0$. 
An exact functor $F:\cA\rightarrow\cB$ between  abelian categories  induces a homomorphism $[F]:K_0(\cA)\rightarrow K_0(\cB)$ of Grothendieck groups.  
Let $A$ be a ring which is free as an abelian 
group,  and $\ba=\{ a_i\}_{i\in I}$  a basis of $A$, such that 
$a_i a_j = \sum_k c_{ij}^k a_k$, where all $c_{ij}^k$ are non-negative
integers. Let further $B$ be an $A$-module.

\begin{definition} \label{def-1} 
A (weak) abelian categorification of $(A,\ba,B)$ consists 
of an abelian category $\mc{B},$  an isomorphism $K_0(\mc{B}) \cong B$ and  
exact endofunctors $F_i: \mc{B} \lra \mc{B}$, $i\in I$, 
such that there are isomorphisms $F_i F_j \cong \oplusop{k} F_k^{c_{ij}^k}$
of functors for all $i,j$,  and for all $i$ the diagram below commutes.
 $$\begin{CD}
    K_0(\mc{B})  @>{[F_i]}>>     K_0(\mc{B})   \\
   @V{\cong}VV                   @VV{\cong}V   \\
    B   @>{a_i}\cdot>>     B
  \end{CD} $$
\end{definition} 

Of course, the existence of a categorification of $B$ implies that there must be a basis of $B$ such that the structure constants are all non-negative 
integers.  
Several known results in the literature can be considered as 
categorifications of modules over various rings, see for example
\cite{BFK}, \cite{CR}, \cite{FKS} and references therein. 

If $R$ is a finite-dimensional algebra over a field,  the Grothendieck group 
$K_0(R\mathrm{-mod})$ of the category of finite-dimensional 
$R$-module is a finite rank free abelian group with the distinguished basis 
given by images of simple $R$-modules. The group $K_0(R\mathrm{-mod})$
has a subgroup $K_0^P(A\mathrm{-mod})$, generated by $[P]$ where
$P$ is a projective $R$-module. $K_0^P(R\mathrm{-mod})$ has the 
distinguished spanning set  given by images of indecomposable projective 
$R$-modules. If the global dimension of $R$ is finite, we have 
$K_0(R\mathrm{-mod})=K_0^P(R\mathrm{-mod})$ but in general
there is only the obvious inclusion
$K_0^P(R\mathrm{-mod})\subset K_0(R\mathrm{-mod}).$ 
For the algebras we consider here this inclusion becomes an isomorphism after 
tensoring with $\Q.$  

Definition~\ref{def-1} also makes sense with $K_0^P$ in place of $K_0,$ 
assuming that the functors $F_i$ take projectives to projectives.


\section{Partitions, Specht modules and Category $\mathcal{O}$} 

Let $n$ be a positive integer  and $S_n$ the symmetric group of order $n!$. 
It is generated by the transpositions $s_i=(i,i+1)$, 
$i=1,\dots,n-1$. Let $\Z[S_n]$ be the integral group algebra of $S_n$. It 
is generated by the element $1$ and the elements 
$\tilde{s}_i=1+s_i$, where $i=1,\dots,n-1$.
 For any {\it partition} $\lambda$ of $n$, i.e.
 $\lambda=(\lambda_1,\dots,\lambda_k)$, $\lambda_1+\dots+\lambda_k=n$, 
$\lambda_i\in\N$, 
$\lambda_1\geq\lambda_2\geq \dots$, we denote by $\mathcal{S}(\lambda)$ the
 right Specht $\mathbb{Z}[S_n]$-module  corresponding to $\lambda$. Note that 
these Specht modules give rise to a complete list of irreducible $S_n$-modules
 over the field $\mathbb{C}$. 
Consider the Lie algebra $\mathfrak{sl}_n(\C)$ with the standard  triangular 
decomposition $\mathfrak{sl}_n(\C)=\mathfrak{n}_-\oplus\mathfrak{h}\oplus 
\mathfrak{n}_+$. Denote by $\mathcal{O}_0$ the corresponding principal block 
of the BGG-category $\mathcal{O}$ for $\mathfrak{sl}_n(\C)$, see \cite{BGG}. 
The simple modules in this category are highest weight modules, with highest 
weights contained in a certain $S_n$-orbit. In particular, the simple objects 
in $\mathcal{O}_0$ are indexed by $w\in S_n$ (we assume that the identity 
element in $S_n$ indexes the one-dimensional  simple module in 
$\mathcal{O}_0$). For $w\in S_n$ we denote by $L(w)$, $\Delta(w)$ and $P(w)$ 
the corresponding simple, Verma and projective modules in $\mathcal{O}_0,$ 
respectively. 

Let $\mu=(\mu_1,\dots,\mu_k)$ be a {\it composition} of $n$, obtained by 
permuting the components of  $\lambda$. Associated with it we have the 
parabolic subalgebra $\mathfrak{p}_\mu$ consisting of $\mu$-block 
upper-triangular matrices in $\mathfrak{sl}_n$, and a subgroup 
$S_\mu=S_{\mu_1}\times S_{\mu_2}\times\cdots\times S_{\mu_n}$ of $S_n$ as 
follows. The algebra $\mathfrak{p}_{\mu}$ contains the Borel algebra 
$\mathfrak{h}\oplus \mathfrak{n}_+$; the Levi factor of $\mathfrak{p}_{\mu}$ 
is the reductive Lie algebra of traceless matrices in  
$\left(\mathfrak{gl}_{\mu_1}(\C)\oplus\dots\oplus 
\mathfrak{gl}_{\mu_k}(\C)\right)+\mathfrak{h}\subset \mf{gl}_n$. The group 
$S_{\mu}$ is the Weyl group of the subalgebra
$\mathfrak{sl}_{\mu_1}(\C)\oplus\dots\oplus \mathfrak{sl}_{\mu_k}(\C)$.

We consider the principal block $\mathcal{O}_0^{\mu}$  of the 
parabolic subcategory  of $\mathcal{O}_0$ associated with 
$\mathfrak{p}_{\mu}$ (see \cite{RC}). By definition this is the full 
subcategory given by all objects $M\in \mathcal{O}_0$ such that 
any $m\in M$ is contained in a {\it finite dimensional} $\p_\mu$-stable 
subspace of $M$ (such modules are called {\em $\p_\mu$-locally finite}). 
The simple objects in 
$\mathcal{O}_0^{\mu}$  are exactly the  simple objects in $\mathcal{O}_0$ 
which correspond to the shortest length representatives in the coset 
$S_{\mu}\setminus S_n$. We denote the latter set by $S(\mu)$ and 
let $w_{\mu}$ denote the longest element in this set. 
For $w\in S(\mu)$ we denote by $L^{\mu}(w)$, $\Delta^{\mu}(w)$ and  
$P^{\mu}(w)$ the  corresponding simple, generalized Verma and 
projective modules in $\mathcal{O}_0^{\mu}$ 
respectively. Remark that $L^{\mu}(w)=L(w)$ for $w\in S(\mu)$.

To a finite dimensional $\mathfrak{sl}_n$-module $E$ we associate the functor
 $F_E:\mathcal{O}_0\rightarrow\mathcal{O}_0$ given by first tensoring with 
$E$ and then taking the largest direct summand contained in $\mathcal{O}_0$. 
In \cite{BG} the indecomposable direct summands of such functors where
 classified (up to isomorphism). They are in natural bijection with the 
elements of $S_n$, where $w\in S_n$ corresponds to the indecomposable 
functor $\theta_w$ which maps $P(e)$ to $P(w)$. For $i=1\dots,n-1$ we 
abbreviate $\theta_i=\theta_{s_i}:\mathcal{O}_0\to
\mathcal{O}_0$. This is the famous {\it translation functor through the 
$s_i$-wall}. All the functors $\theta_w$ are exact and, by definition, 
preserve $\mathcal{O}_0^{\mu}$. The translations $\theta_i$ through walls are 
self-adjoint functors. 

For $w\in S_n$ let $\mathcal{R}(w)$ be the right cell of the 
element $w$ (see for example \cite{KL}, or \cite{Douglassinv} for a 
characterization of $\mathcal{R}(w)$). According to the main result of 
\cite{Ir1}, the modules $P^{\mu}(w)$, $w\in \mathcal{R}(w_{\mu})$, constitute 
an 
exhaustive list of indecomposable projective-injective modules in 
$\mathcal{O}_0^{\mu}$. 

\section{A categorification of the Specht modules}

Set $P_{\mu}=\oplus_{w\in \mathcal{R}(w_{\mu})}P^{\mu}(w)$ and denote 
by $\mc{C}_{\mu}$ the full subcategory of $\mathcal{O}_0^{\mu}$ 
which consists of all $M$ admitting a two-step resolution, 
\begin{equation}\label{eq1}
P_1\to P_0\to M\to 0, 
\end{equation}
with $P_1,P_0\in\mathrm{add}(P_{\mu})$. By \cite[Section~5]{Au}, 
the category $\mc{C}_{\mu}$ is equivalent to $A_{\mu}\mathrm{-mod}$, where 
$A_{\mu}=\mathrm{End}_{\mathcal{O}_0^{\mu}}\left(P_{\mu}\right)$, in particular,
$\mc{C}_{\mu}$ is an abelian category. Since projective functors preserve both, 
projectivity and injectivity, they preserve $\mathrm{add}(P_{\mu})$. Now 
from \eqref{eq1} it follows that they preserve $\mc{C}_{\mu}$ as well. 

\begin{theorem}\label{t1}
The action of $\theta_i$, $i=1,\dots,n-1$, on $\mc{C}_{\mu}$ 
defines on $K_0(\mc{C}_{\mu})$ the structure of an $S_n$-module.
The restriction of this action to 
$K_0^P(\mc{C}_{\mu})$ categorifies $\mathcal{S}(\lambda')$. 
The Kazhdan-Lusztig basis in $\mathcal{S}(\lambda')$
corresponds to the basis given by the isomorphism classes of the 
indecomposable projective modules, and the elements
$\theta_i=\theta_{s_i}$ correspond to the generators $\tilde{s}_i$ of
$\mathbb{Z}\, S_n$.
\end{theorem}

The associative algebras which correspond to both $\mathcal{O}_0$ and
$\mathcal{O}_0^{\mu}$ are Koszul (see \cite{BGS}) 
and admit a canonical
positive grading (the Koszul grading), which we fix. This allows us
to consider graded versions of both $\mathcal{O}_0$ and
$\mathcal{O}_0^{\mu}$, see \cite{St} for details. If some confusion 
can arise, to emphasize the graded versions, we will add the superscript 
$gr.$ Let $A_\mu-\operatorname{gmod}$ denote the category of all finitely 
generated graded $A_\mu$-modules, where $A_\mu$ is equipped with the 
grading 
induced from the Koszul grading mentioned above.  We denote by
$\langle 1\rangle:A_\mu-\operatorname{gmod}\rightarrow 
A_\mu-\operatorname{gmod}$ the functor which increases the grading by 
$1$. 

In \cite{St} it was
shown that simple modules, Verma modules, generalized Verma modules and 
projective modules in
both $\mathcal{O}_0$ and $\mathcal{O}_0^{\mu}$ are gradable and, moreover,
that the functors $\theta_i$ are gradable as well. Let
 $\theta_i:A_\mu-\operatorname{gmod}\rightarrow A_\mu-\operatorname{gmod}$ be 
the standard graded lift of the translation through the wall (see \cite{St}). 

Let $\mathcal{H}_n$ denote the Hecke algebra of $S_n$ over 
$\mathbb{Z}[v,v^{-1}]$ with the $\Z[v,v^{-1}]$-basis  $H_{w}$, $w\in W$, and 
the relations $H_s^2=H_e+(v-v^{-1})H_s$, $H_sH_t=H_tH_s$ if $st=ts$, and 
$H_sH_tH_s=H_tH_sH_t$ if $sts=tst$, where $s, t\in\{s_i\mid 1\leq i<n\}$. 
The algebra $\mathcal{H}_n$ is a deformation of $\Z[S_n]$.  
Let $\tilde{\mathcal{S}}(\lambda)$ be the (right) Specht 
$\mathcal{H}_n$-module, which corresponds to $\lambda$ (see e.g. 
\cite{Murphy} or \cite{KL}). 

\begin{theorem}\label{t2}
The graded functors $\theta_i$, $i=1,\dots,n-1$ and the functor $\langle 1\rangle$ induce endofunctors of $\mc{C}_{\mu}^{gr}$ and define on $K_0(\mc{C}_{\mu}^{gr})$ the  
structure of an $S_n$-graph. The restriction of this
action to $K_0^P(\mc{C}_{\mu}^{gr})$ categorifies
$\tilde{\mathcal{S}}(\lambda'),$ and the 
Kazhdan-Lusztig basis of $\tilde{\mathcal{S}}(\lambda')$ corresponds to the basis 
given by the iso-classes of the indecomposable  projective modules 
with head concentrated in degree zero, the action of 
$\theta_i$ (resp. $\langle 1\rangle$) corresponds to the action of 
$H_{s_i}+v$ (resp. of $v$).
\end{theorem}

\begin{proof}[Proof of Theorem~\ref{t1} and Theorem~\ref{t2}.]
Theorem~\ref{t1} is obtained from Theorem~\ref{t2} via evaluation $v=1$,
hence it is enough to prove Theorem~\ref{t2}. As a consequence of the
Kazhdan-Lusztig Theorem and \cite[Theorem~5.1]{St}, for $w\in S_n$ we 
have the following equality in $K_0(\mathcal{O}_0^{gr})$:
\begin{equation}\label{eq5}
[\theta_i L(w)]=
\begin{cases}
0 & ws_i>w,\\
[L(w)\langle 1\rangle]+[L(w)\langle -1\rangle]
+\sum_{w'\in S_n}\mu(w',w) [L(w')] & ws_i<w,
\end{cases}
\end{equation}
where $\mu(w',w)$ denotes the Kazhdan-Lusztig $\mu$-function, see 
\cite[Definition~1.2]{KL}. 
If $w\in \mathcal{R}(w_{\mu})$ then we get the following equality for the
induced action on  $K_0(\mc{C}_{\mu}^{gr})$:
\begin{equation}\label{eq6}
[\theta_i L^{\mu}(w)]=
\begin{cases}
0 & ws_i>w,\\
[L^{\mu}(w)\langle 1\rangle]+[L^{\mu}(w)\langle -1\rangle]
+\sum_{w'\in \mathcal{R}(w_{\mu})}\mu(w',w) [L^{\mu}(w')] & ws_i<w.
\end{cases}
\end{equation}

This means that the action of $\theta_i$  and $\langle -1\rangle$ on 
$\mc{C}_{\mu}^{gr}$ defines on $G_0(\mc{C}_{\mu}^{gr})$ the structure of an $S_n$-graph. By \cite[Theorem~1.4]{KL}, we obtain a simple $\mathcal{H}_n$-module
over any field of characteristic $0$. So, to complete the proof we just 
have to show that this module is exactly $\tilde{\mathcal{S}}(\lambda')$.
And to do this we can specialize to $v=1$. In this case the necessary
statement follows from \cite[Theorem~4.1]{Na}. 
\end{proof}

\begin{Example}\label{ex1}
{\rm {\bf The trivial representation.}
If $\mu=(1,1,1,\dots,1)$ then $\mathcal{O}_{0}^{\mu}=\mathcal{O}_{0}$. 
There is a unique indecomposable projective-injective module $P(w_0)$, 
where $w_0$ is the longest element in $S_n$.  We have
$\theta_i P(w_0)\cong P(w_0)\langle 1\rangle\oplus P(w_0)\langle-1\rangle$, which implies
$[\theta_i][L(w_0)]=[L(w_0)\langle 1\rangle]+[L(w_0)\langle-1\rangle]$. Hence
$K_0(\mc{C}^{gr}_{\mu})\cong\mathbb{Z}[v,v^{-1}]$ and also $K^P_0(\mc{C}^{gr}_{\mu})\cong\mathbb{Z}[v,v^{-1}]$
(note that $K^P_0(\mc{C}^{gr}_{\mu})\subsetneq K_0(\mc{C}^{gr}_{\mu})$
in this case); the action of 
$H_{i}+v$ is just the multiplication by $v+v^{-1}$. 
This construction therefore categorifies the trivial representation 
of the symmetric  group and its Hecke algebra deformation.

We remark that in this case the endomorphism algebra
$A_{\mu}$ of $P(w_0)$ is isomorphic to the complex cohomology ring of 
the full flag variety of $\C^n,$ and $\cA$ is equivalent to the category 
of graded finite-dimensional modules over this algebra.
Further, $\theta_i$ is  isomorphic to the functor of tensoring 
with the graded $A_{\mu}$-bimodule which is the cohomology ring of an 
iterated flag variety (having  two subspaces of dimension $i$). 
}
\end{Example}

\begin{Example}\label{ex2}
{\rm {\bf The sign representation.}
If $\mu=(n)$ then $\mathcal{O}_{0}^{\mu}$ is semi-simple and
contains a unique, up to isomorphism, simple projective (and hence 
also injective) module, namely the trivial module $L(e)$. We have
$\theta_i L(e)=0,$ so that $K_0(\mc{C}_{\mu})=K_0^P(\mc{C}_{\mu})\cong 
\mathbb{Z}[v,v^{-1}]$ and $H_i+v$ acts by $0.$ We obtain a categorification of the 
sign representation of the symmetric group and its Hecke algebra 
deformation.}
\end{Example}

\begin{Remark}\label{remark1}
{\rm 
The action of projective functors on (the graded version of) $\mathcal{O}_0$
categorifies the regular representation of $S_n$ (respectively 
$\mathcal{H}_n$). This is a straightforward calculation using the basis
$[\Delta(w)]$, $w\in S_n$, of $K_0(\mathcal{O}_0)$. The
basis $[P(w)]$, $w\in S_n$, then corresponds to the
Kazhdan-Lusztig basis in $\mathbb{Z}[S_n]$ (respectively 
$\mathcal{H}_n$).
}
\end{Remark}

\section{The bilinear forms} 

Now recall that every Specht module $\tilde{\mathcal{S}}(\lambda)$ has a 
symmetric, non-degenerate, $\mathcal{H}_n$-invariant bilinear form $<,>$ with 
values in $\Z[v,v^{-1}]$ which is unique up to a scalar 
\cite[page~114]{Murphy}. We give a categorical interpretation of this form. 
For any $\Z$-graded complex vector space $M=\oplus_{j\in\Z} M^j$ let 
$h(M)=\sum_{j\in\Z}
(\operatorname{dim}_\mC M^j)v^{j}\in\Z[v,v^{-1}]$ be the corresponding 
Hilbert polynomial. 

Consider $A_\mu-\operatorname{gmod}$, the graded version of $\mc{C}_{\mu}$. 
If $M$, $N\in A_\mu-\operatorname{gmod}$ then the space $E^i(M,N):=
\operatorname{Ext}^i_{A_\mu-\operatorname{mod}}(M,N)$ is a $\Z$-graded vector 
space for any $i\in\Z$. Set $h(E(M,N))=\sum_{i\in\Z}(-1)^i h(E^i(M,N))$. 
Let $\operatorname{d}$ denote the 
graded lift of the standard duality on $\cO_0^\mu$, restricted to the 
category $\mc{C}^{gr}_{\mu}$ (see \cite{St}).

The following result categorifies the bilinear forms on Specht modules.

\begin{prop}
Let $\mu$ be a composition of $n$ and $\la$ be the corresponding partition. 
Then $$\beta(_-, _-):=
h(E(_- ,\operatorname{d}(_-))):\mathrm{Ob}(A_\mu-\operatorname{gmod})\times 
\mathrm{Ob}(A_\mu-\operatorname{gmod})\rightarrow \Q(v)$$ descends to a symmetric, 
non-degenerate, $\mathcal{H}_n$-invariant bilinear form $<\cdot,\cdot>$ on 
$K_0(A_\mu-\operatorname{gmod})$ and 
$K_0^P(A_\mu-\operatorname{gmod})$ with values in $\Q(v)$ and
$\Z[v,v^{-1}]$ respectively.
\end{prop}

\begin{proof}
The form is bilinear, since the bifunctor $\beta$ is additive in both 
arguments. To verify the symmetry, it is therefore enough to check that 
$\beta(M,N)=\beta(N,M)$, when $M$, $N$ are projective (and therefore 
injective). We have $\beta(M,N)=\beta(\operatorname{d}(N), 
\operatorname{d}(M))=\beta(N,M)$, since 
$\operatorname{d}(M)\cong M\langle k\rangle$ and 
$\operatorname{d}(N)\cong N\langle k\rangle$ for some $k\in\Z$,
which is the same for both $M$ and $N$ by \cite[Theorem~5.2(1)]{MS}. So, the 
form is symmetric. 
To see that the form is $\mathcal{H}_n$-invariant we only have to verify that 
$\beta(P^\mu(x), \theta_i P^\mu(w))=(v+v^{-1}) \beta(P^\mu(x), P^\mu(w))$ for
 any $x, w\in S_n$ such that $xs_i<x$ and $ws_i<w$ in the Bruhat ordering 
(\cite[Proposition~6.1]{Fung}). This is however clear, since 
$ \theta_i P^\mu(w)\cong P^\mu(w)\langle 1\rangle \oplus  
P^\mu(w)\langle -1\rangle$ (this follows from \cite[Theorem~5.1]{St}). 

To see that it is non-degenerate (i.e. that it has trivial radical), we 
observe that in the basis given by indecomposable projectives with heads 
concentrated in degree zero, the form is given by the matrix 
$M=h(\operatorname{Hom}_{A_\mu}(_-,\operatorname{d}(_-)))$. Since the algebra 
$A_\mu$ is positively graded, it follows that the entries off the diagonal of 
$M$ are in $v\Z[v]$, and the ones on the diagonal are of the form $1+f$, where
 $f\in v\Z[v]$. Hence, the matrix $M$ has non-zero determinant and thus
the form is non-degenerate.     

It is left to note that the algebra $A_{\mu}$ has infinite 
global dimension in general, and hence $h(E(M,N))$ is a  power series and 
not a Laurent polynomial for general $M$ and $N$. However, since the
(graded) Cartan matrix of $A_{\mu}$ is invertible, it follows that
for simple $M$ and $N$ we can compute $h(E(M,N))$ by inverting the
Cartan matrix and hence $h(E(M,N))$ can be written as a rational 
function. 
\end{proof}

\end{document}